\documentclass[12pt,twoside]{article}
\title{On Weakly Coherent Rings}
\date{}
\usepackage{amsfonts}
\usepackage{amsmath}
\usepackage{amssymb}
\usepackage{latexsym}
\newtheorem{thm}{\bf Theorem}[section]

\newtheorem{lem}[thm]{\bf Lemma}
\newtheorem{prop}[thm]{\bf Proposition}

\newtheorem{rem}[thm]{\bf Remark}

\newtheorem{exmp}[thm]{\bf Example}

\catcode`\ç=13
\defç{\c{c}}
\catcode`\é=13
\defé{\'e}
\catcode`\à=13
\defà{\`a}
\catcode`\è=13
\defè{\`e}
\catcode`\â=13
\defâ{\^a}
\catcode`\ù=13
\defù{\`u}
\catcode`\ê=13
\defê{\^e}
\catcode`\î=13
\defî{\^\i}
\catcode`\ô=13
\defô{\^o}
\newcommand{\field}[1]{\mathbb{#1}}

\newcommand{\R}{\field{R}}
\newcommand{\Q }{\field{Q}}
\newcommand{\Z }{\field{Z}}
\newcommand{\N }{\field{N}}
\def\proof{{\parindent0pt {\bf Proof.\ }}}

\begin{document}
\thispagestyle{empty}
%%%%%%%%%%%%%%%%%%%%%%%%%%%%%%%%%%%%%%%%%%%%%%%%%%%%%%%%%
%%%%%%%%%%%%%%%%%%%%%%%%%%%%%%%%%%%%%%%%%%%%%%%%%%%%%%%%%
%%%%%%%%%%%%%%%%%%%%%%%%%%%%%%%%%%%%%%%%%%%%%%%%%%%%%%%%%
%%%TITLE%%%%%%%%%%%%%%%%%%%%%%%%%%%%%%%%%%%%%%%%%%%%%%%%%
\maketitle \vspace*{-2cm}
\begin{center}{\large\bf Chahrazade Bakkari and Najib Mahdou}
%%%%%%%%%%%%%%%%%%%%%%%%%%%%%%%%%%%%%%%%%%%%%%%%%%%%%%%%%
%%%%%%%%%%%%%%%%%%%%%%%%%%%%%%%%%%%%%%%%%%%%%%%%%%%%%%%%%
%%%%%%%%%%%%%%%%%%%%%%%%%%%%%%%%%%%%%%%%%%%%%%%%%%%%%%%%%
%%%NAMES%%%%%%%%%%%%%%%%%%%%%%%%%%%%%%%%%%%%%%%%%%%%%%%%%

\bigskip
%%%%%%%%%%%%%%%%%%%%%%%%%%%%%%%%%%%%%%%%%%%%%%%%%%%%%%%%%
%%%%%%%%%%%%%%%%%%%%%%%%%%%%%%%%%%%%%%%%%%%%%%%%%%%%%%%%%
%%%%%%%%%%%%%%%%%%%%%%%%%%%%%%%%%%%%%%%%%%%%%%%%%%%%%%%%%
%%%%%%%%%%%%ADDRESSES%%%%%%%%%%%%%%%%%%%%%%%%%%%%%%%%%%%%%%%%%%%%%
\small{Department of Mathematics, Faculty of Science and
Technology of Fez,\\ Box 2202, University S. M. Ben Abdellah Fez,
Morocco \\ cbakkari@hotmail.com \\ mahdou@hotmail.com }
\end{center}

\bigskip\bigskip
%%%%%%%%%%%%%%%%%%%%%%%%%%%%%%%%%%%%%%%%%%%%%%%%%%%%%%%%%
%%%%%%%%%%%%%%%%%%%%%%%%%%%%%%%%%%%%%%%%%%%%%%%%%%%%%%%%%
%%%%%%%%%%%%%%%%%%%%%%%%%%%%%%%%%%%%%%%%%%%%%%%%%%%%%%%%%
%%%ABSTRACT%%%%%%%%%%%%%%%%%%%%%%%%%%%%%%%%%%%%%%%%%%%%%%
\noindent{\large\bf Abstract.}   In this paper, we define weakly
coherent rings, and examine the transfer of these rings to
homomorphic image, trivial ring extension, localization, and
direct product. These results provide examples of weakly coherent
rings that are not coherent rings. We show that the class of
weakly coherent rings is not stable by localization. Also, we show
that the class of weakly coherent rings and the class of strongly
$2$-coherent rings are not comparable.
\bigskip

%%%%%%%%%%%%%%%%%%%%%%%%%%%%%%%%%%%%%%%%%%%%%%%%%%%%%%%%%
\small{\noindent{\bf Key Words.}  Weakly coherent ring, coherent
ring, homomorphic image, trivial ring extension, localization, direct product. \\
\bigskip

%%%%%%%%%%%%%%%%%%%%%%%%%%%%%%%%%%%%%%%%%%%%%%%%%%%%%%%%%
%%%%%%%%%%%%%%%%%%%%%%%%%%%%%%%%%%%%%%%%%%%%%%%%%%%%%%%%%
%%%%%%%%%%%%%%%%%%%%%%%%%%%%%%%%%%%%%%%%%%%%%%%%%%%%%%%%%
\small{\noindent{\bf 2000 Mathematics Subject Classification: 13F05, 13B05, 13A15, 13D05, 13B25.}    \\
%%%%%%%%%%%%%%%%%%%%%%%%%%%%%%%%%%%%%%%%%%%%%%%%%%%%%%%%%
%%%INTRODUCTION%%%%%%%%%%%%%%%%%%%%%%%%%%%%%%%%%%%%%%%%%%
\begin{section}{Introduction}

Throughout this paper all rings
are assumed to be commutative with identity elements and all
modules are unital.\\
Let $R$ be a commutative ring. We say that an ideal is regular if it contains a regular
element, i.e; a non-zerodivisor element.   \\

For a nonnegative integer $n$, an $R$-module $E$ is $n$-presented
if there is an exact sequence of $R$-modules: \\
$$F_n \rightarrow  F_{n-1} \rightarrow \ldots F_1  \rightarrow F_0  \rightarrow E
 \rightarrow 0$$
 where  each $F_i$ is a finitely generated free $R$-module. In
particular, $0$-presented and $1$-presented $R$-modules are,
respectively, finitely generated and finitely presented
$R$-modules. \\
A ring $R$ is coherent if every finitely generated ideal of $R$ is
finitely presented; equivalently, if $(0:a)$ and $I\cap J$ are
finitely generated for every $a\in R$ and any two finitely
generated ideals $I$ and $J$ of $R$ \cite{Gz2}. Examples of
coherent rings are Noetherian rings, Boolean algebras, von Neumann
regular rings, valuation rings, and Pr\"ufer/semihereditary rings. See for instance \cite{Gz2}.\\

In this paper, we investigate a particular class of coherent rings
that we call weakly coherent rings. A ring $R$ is called a weakly
coherent ring if for each proper ideals $I \subseteq J$ such that
$I$ is finitely generated and $J$ is finitely presented, then $I$
is finitely presented. If $R$ is coherent, then $R$ is naturally weak coherent.
Our aim in this paper is to prove that the converse is false in general. \\

  We say that $R$ is strong $n$-coherent if each
$n$-presented $R$-module is $(n+1)$-presented. In particular, any
coherent ring (i.e., $1$-coherent ring) is a strong $2$-coherent
ring. This led us to consider the relation between the class of
weakly coherent rings and the class of strong $2$-coherent rings.\\

Let $A$ be a ring, $E$ be an $A$-module and $R :=A \propto E$ be
the set of pairs $(a,e)$ with pairwise addition and multiplication
given by $(a,e)(b,f) =(ab,af+be)$. $R$ is called the trivial ring
extension of $A$ by $E$. Considerable work has been concerned with
trivial ring extensions. Part of it has been summarized in Glaz's
book \cite{Gz2}) and Huckaba's book (where $R$ is called the
idealization of $E$ by $A$) \cite{H}). \\

In the context of non-total ring of quotients (i.e., ring
containing a regular element), we show that the notion of weak
coherent coincide with the definition of coherent ring. The goal
of this work is to exhibit a class of non-coherent weakly coherent
rings. We show that the class of  weakly coherent rings is not
stable by localization. Also, we show that the class of  weakly
coherent rings and the class of strong $2$-coherent rings are not
comparable. For this purpose, we study the transfer of this
property to homomorphic image, trivial ring extension, and direct product. \\

\end{section}
%%%%%%%%%%%%%%%%%%%%%%%%%%%%%%%%%%%%%%%%%%%%%%%%%%%%%%%%%%%%
%%%%%%%%%%%%%%%%%%%%%%%%%%%%%%%%%%%%%%%%%%%%%%%%%%%%%%%%%
%%%%%%%%%%%%%%%%%%%%%%%%%%%%%%%%%%%%%%%%%%%%%%%%%%%%%%%%%
%%%%%%%%%%%%%%%%%%%%%%%%%%
%%%%%%%%%%%%%%%%%%%%%        Section 2: Main Results
%%%%%%%%%%%%%%%%%%%%%
%%%%%%%%%%%%%%%%%%%%%%%%%%%%%%%%%%%%%%%%%%%%%%%%%%%%%%%%%%%%%
%%%%%%%%%%%%%%%%%%%%%%%%%%%%%%%%%%%%%%%%%%%%%%%%%%%%%%%%%%%%%%%
\begin{section}{Main Results}

Recall that for nonnegative integers $n$ and $d$, we say that a
ring $R$ is an $(n, d)$-ring if $pd_{R}(E) \leq d$ for each
$n$-presented $R$-module $E$ (as usual,  $pd_{R} E$ denotes the
projective dimension of $E$ as an $R$-module). See for instance \cite{C,KM1,KM2,M1,M2}.\\
We begin this section by giving an example of non-coherent weakly coherent ring.\\

\bigskip

%%%%%%%%%%%%%%%%%%%%%%%%%%%%%%%%%%%%%%%%%%%%%%%%%%%%%%%%%%%%%%%%%%%%%%%%%%%%%%%
%%%%%%%%%%%%%%%%%%%%%%%%%%%%%%%%%%%%%%%%%%%%%%%%%%%%%%%%%%%%%%%%%%%%%%%%%%%%%%%
\begin{exmp}\label{result1}
 Let $(A,M)$ be a local ring with maximal ideal $M$, $E$ be an
 $A/M$-vector space with infinite rank, and let $R :=A \propto E$  be the trivial ring extension of $A$
by $E$.  Then: \\
{\bf 1)} $R$ is a weak coherent ring.\\
{\bf 2)} $R$ is not a coherent ring. $\Box$
\end{exmp}

\proof {\bf 1)} It suffices to show that there is no finitely
presented proper ideal $J$ of $R$. Deny. Let $J$ be a finitely
presented proper ideal of $R$. Then $J$ is free since $R$ is a
local $(2,0)$-ring by \cite[Theorem 2.1(1)]{M2}, that is $J =Ra$
for some regular element $A$ of $R$; a contradiction since $J
\subseteq M \propto E$ and
$(M \propto E)(0,e) =(0,0)$ for each $e \in E-\{0\}$. Hence $R$ is a weak coherent ring. \\
{\bf 2)} We claim that $R$ is not coherent. Deny. Assume that $R$ is coherent. But, $R$
is a $(2,0)$-ring by \cite[Theorem 2.1(1)]{M2}. Hence, $R$ is Von Neumann regular ring since it is coherent, a contradiction
by \cite[Theorem 2.1(2)]{M2}. Hence, $R$ is not coherent, as desired. $\Box$ \\

\bigskip

%%%%%%%%%%%%%%%%%%%%%%%%%%%%%%%%%%%%%%%%%%%%%%%%%%%%%%%%%%%%%%%%%%%%%%%%%%%%%%
%%%%%%%%%%%%%%%%%%%%%%%%%%%%%%%%%%%%%%%%%%%%%%%%%%%%%%%%%%%%%%%%%%%%%%%%%%%%%%%

Now, we give a sufficient condition to have equivalence between a coherence and weakly coherence properties. \\

\bigskip

\begin{prop}\label{result2}
 Let $R$ be a commutative ring, then: \\
{\bf 1)} If $R$ is a coherent ring, then $R$ is a weak coherent ring.\\
{\bf 2)} Assume that $R$ contains a regular element (that is $R$
is not a total ring of quotients). Then $R$ is a coherent ring if
and only if $R$ is a weak coherent ring. $\Box$
\end{prop}

\proof {\bf 1)} Clear. \\
{\bf 2)} It remains to show that if $R$ is a weak coherent ring
and contains a regular element $a$, then $R$ is a coherent ring.
Let $I $ be a finitely generated proper ideal of $R$. Hence, $aI
\subseteq aR$, $aI$ is a finitely generated proper ideal of $R$,
and $aR (\cong R)$ is a finitely presented proper ideal of $R$.
Therefore, $aI$ is a finitely presented ideal of $R$ and so $I
(\cong aI$) (since $a$ is regular) is finitely presented, as
desired. $\Box$

\bigskip

\begin{rem}\label{result3}
By the above result, a non-coherent weak coherent ring is
necessary a total ring of quotients. $\Box$
\end{rem}

%%%%%%%%%%%%%%%%%%%%%%%%%%%%%%%%%%%%%%%%%%%%%%%%%%%%%%%%%%%%%%%%%%%%%%%%%%%%%%
%%%%%%%%%%%%%%%%%%%%%%%%%%%%%%%%%%%%%%%%%%%%%%%%%%%%%%%%%%%%%%%%%%%%%%%%%%%%%%%

\bigskip

Now, we investigate the homomorphic image of weak coherent
rings.\\

\bigskip

 \begin{thm}\label{result4}
 Let $R$ be a weak coherent ring and $I$ be a finitely generated
 ideal of $R$. Then $R/I$ is a weak coherent ring.
\end{thm}

\proof Let $L/I \subseteq J/I$ be two finitely generated proper
ideals of $R/I$ such that $J/I$ is finitely presented. Our aim is
to show that $L/I$ is finitely presented. \\

Remark that $L \subseteq J$ are two finitely generated proper
ideals of $R$. We claim that $J$ is finitely presented. \\
Indeed, there exists an exacte sequence of $(R/I$)-modules:
$$0 \rightarrow T \rightarrow (R/I)^n \rightarrow  J/I
\rightarrow 0 (*)$$ where $n$ is a positive integer and $T$ is a
finitely generated $(R/I)$-module (since $J/I$ is a finitely
presented ideal of $R/I$). Hence, $T$ is a finitely generated
$R$-module. On the other hand, $R/I$ is a finitely presented
$R$-module (since $I$ is a finitely generated ideal of $R$ and by
the exact sequence of $R$-modules $0 \rightarrow I \rightarrow R
\rightarrow  R/I \rightarrow 0$). Therefore, $J/I$ is a finitely
presented $R$-module by exact sequence $(*)$ considered as an
exact sequence of $R$-modules. Consequently, $J$ is a finitely
presented ideal of $R$ by the exact sequence of $R$-modules $0
\rightarrow I \rightarrow J \rightarrow  J/I \rightarrow 0$, as
desired. \\

Now, we have $L \subseteq J$, $L$ is a finitely generated ideal,
and $J$ is a finitely presented ideal; so $L$ is a finitely
presented ideal of $R$ since $R$ is a weakly coherent ring. Hence,
the exact sequence of $R$-modules $0 \rightarrow I \rightarrow L
\rightarrow L/I \rightarrow 0$ shows that $L/I$ is a finitely
presented $R$-module. We claim that $L/I$ is a finitely presented
ideal of $R/I$ and this completes the proof of Theorem~\ref{result4}. \\
Indeed, since $L/I$ is a finitely generated ideal of $R/I$,
consider the exact sequence of $(R/I)$-modules: \\
$$0 \rightarrow S \rightarrow (R/I)^m \rightarrow  L/I \rightarrow 0 (**)$$
where $m$ is a positive integer and $S$ is an $(R/I)$-module. The
exact sequence $(**)$ is also an exact sequence of $R$-modules;
hence $S$ is a finitely generated $R$-module since $L/I$ and $R/I$
are finitely presented $R$-modules. Therefore, $S$ is also a
finitely generated $(R/I)$-module and the exact sequence of
$(R/I)$-modules $(**)$ shows that $L/I$ is a finitely presented
ideal of $R/I$, as desired. $\Box$

\bigskip

The condition "$I$ is a finitely generated ideal of $R$" is
necessary in Theorem~\ref{result4} as the following example shows:

\bigskip

\begin{exmp}\label{result5}
 Let $(A,M)$ be a non-coherent local domain with maximal ideal $M$, $E$ be an
 $A/M$-vector space with infinite rank, let $R :=A \propto E$  be the trivial ring extension of $A$
by $E$, and set $I :=0 \propto E$.  Then: \\
{\bf 1)} $R$ is a weak coherent ring by Example~\ref{result1}.\\
{\bf 2)} $R/I (\cong A)$ is not a weak coherent ring (by
Proposition~\ref{result2} since $A$ is not a coherent domain).
$\Box$
\end{exmp}

%%%%%%%%%%%%%%%%%%%%%%%%%%%%%%%%%%%%%%%%%%%%%%%%%%%%%%%%%%%%%%%%%%%%%%%%%%%%%%
%%%%%%%%%%%%%%%%%%%%%%%%%%%%%%%%%%%%%%%%%%%%%%%%%%%%%%%%%%%%%%%%%%%%%%%%%%%%%%%
%%%%%%%%%%%%%%%%%%%%%%%%%%%%%%%%%%%%%%%%%%%%%%%%%%%%%%%%%%%%%%%%%%%%%%%%%%%%%%
%%%%%%%%%%%%%%%%%%%%%%%%%%%%%%%%%%%%%%%%%%%%%%%%%%%%%%%%%%%%%%%%%%%%%%%%%%%%%%%

\bigskip

Now, we investigate a weak coherent property in a particular class
of total rings of quotients; namely, those arising as trivial ring
extensions of local rings by vector spaces over the residue
fields. The following main result enriches the literature with
original examples of non-coherent weak coherent rings.

 \bigskip

\begin{thm}\label{result6}
 Let $(A,M)$ be a local ring with maximal ideal $M$, $E$ be an
 $A$-module such that $ME =0$, and let $R :=A \propto E$  be the trivial ring extension of $A$
by $E$.  Then $R$ is a weak coherent ring if and only if one of the following two properties holds: \\
{\bf 1)} $E$ is an $A/M$-vector space with infinite rank.\\
{\bf 2)} $E$ is an $A/M$-vector space with finite rank and $A$ is
weak coherent. $\Box$
\end{thm}

\proof Assume that $R$ is a weak coherent ring and $E$ is an
$A/M$-vector space with finite rank. Our aim is to show that $A$
is weak coherent. \\

Let $I \subseteq J$ be two proper ideals of $A$ such that $I$ is
finitely generated and $J$ is finitely presented. Then $(I \propto
0) \subseteq (J \propto 0)$ are two finitely generated proper
ideals of $R$. We claim that $J \propto 0$ is a finitely presented
ideal of $R$. \\
Indeed, let $J :=\sum_{i=1}^{n}Aa_{i}$ for some positive integer
$n$ and some $a_i \in J$, and consider the exact
sequence of $A$-modules: \\
$$0 \rightarrow Ker(u) \rightarrow A^n  \buildrel u \over\rightarrow  J \rightarrow 0$$
where $u$ is defined by $u((b_{i})_{1\leq i\leq n})
=\sum_{i=1}^{n}b_{i}a_{i}$.  Hence, $Ker(u) (:\{(b_{i})_{1\leq
i\leq n} \in A^n / \sum_{i=1}^{n}a_{i}b_{i} =0\}$) is a finitely
generated $A$-module (since $J$ is a finitely presented ideal of
$A$). Now, consider the exact sequence of $R$-modules: \\
$$0 \rightarrow  Ker(v) \rightarrow  R^n  \buildrel v \over\rightarrow
J \propto 0 \rightarrow 0$$ where $v((b_{i},e_{i})_{1\leq i\leq
n}) =\sum_{i=1}^{n}(b_{i},e_{i})(a_{i}, 0)$. But,
\begin{eqnarray}
 \nonumber  Ker(v)&=& \{(b_{i},e_{i})_{1\leq i\leq n} \in R^n / \sum_{i=1}^{n}(a_{i},0)(b_{i},e_{i}) =0\} \\
 \nonumber          &=& \{(b_{i},e_{i})_{1\leq i\leq n} \in R^n / \sum_{i=1}^{n}a_{i}b_{i} =0\} \\
\nonumber         &=& Ker(u) \propto E^{n}
\end{eqnarray}
(since $a_{i} \in J \subseteq M$) which is finitely generated
$R$-module (since $Ker(u)$ and $E$ are finitely generated
$A$-modules). Hence, $J \propto 0$ is a finitely presented
(proper) ideal of $R$ and so $(I \propto 0) (\subseteq (J \propto
0))$ is a finitely presented ideal of $R$ since $R$ is a weak
coherent ring. Therefore, by the same reasoning as for $J$ above,
we can show that $I$ is finitely
presented and this shows that $A$ is a weak coherent ring. \\

Conversely, if $dim_{A/M}E =\infty $, then $R$ is a weak coherent
ring by Example~\ref{result1}. Now, assume that $dim_{A/M}E <
\infty $ and $A$ is a weak coherent ring and our aim is to show
that
$R$ is a weak coherent ring. \\
Let $I (:=\sum_{i=1}^{n}R(a_{i},e_{i})) \subseteq J
(:=\sum_{i=1}^{m}R(b_{i},f_{i}))$ be two proper ideals of $R$ such
that $n, m$ are positive integers, $a_{i}, b_{j} \in A$ and
$e_{i}, f_{j} \in E$ for each $i, j$, and $J$ is finitely
presented. We which to show that $I$ is finitely presented. Two
cases are then possible: \\
{\bf Case 1.} $b_{i} =0$ for all $i =1, \ldots , m$.\\
In this case, $a_{i} =0$ for all $i =1, \ldots , n$ and $I :=0
\propto E_{1}$ and $J :=0 \propto E_{2}$ for some $(A/M)$-vector
subspace $E_1$ and $E_2$ of $E$. Assume that $(e_{i})_{i =1,
\ldots , n}$ and $(f_{i})_{i =1, \ldots , m}$ are respectively
basis of the $(A/M)$-vector space $E_1$ and $E_2$.
Consider the exact sequence of $R$-modules: \\
$$0 \rightarrow  Ker(u) \rightarrow  R^m  \buildrel u \over\rightarrow
J (:= 0 \propto E_{2}) \rightarrow 0$$ where
$u((c_{i},g_{i})_{1\leq i\leq m}) =\sum_{i=1}^{m}(c_{i},g_{i})(0,
f_{i}) =(0, \sum_{i=1}^{m}c_{i}f_{i})$. Hence, $Ker(u) =M^{m}
\propto E^{m} (=(M \propto E)^{m})$ since $(f_{i})_{i =1, \ldots ,
m}$ is a basis of the $(A/M)$-vector space $E_2$ and so $M$ is a
finitely generated ideal of $A$ (since $J$ is a finitely presented
ideal of $R$). Therefore, the exact sequence of $R$-modules: \\
$$0 \rightarrow  M^n \propto E^n \rightarrow  R^n  \buildrel v \over\rightarrow
I (:= 0 \propto E_{1}) \rightarrow 0$$ where
$v((c_{i},g_{i})_{1\leq i\leq n}) =\sum_{i=1}^{n}(c_{i},g_{i})(0,
e_{i})$ shows that $I$ is a finitely presented ideal of $R$ (since
$M^n \propto E^n$ is a finitely generated $R$-module), as desired.
\\
{\bf Case 2.} $b_{i} \not= 0$ for some $i =1, \ldots ,m$. \\
We may assume that $((a_{i},e_{i})_{i =1, \ldots ,n})$ and
$((b_{i},f_{i})_{i =1, \ldots ,n})$ are minimal generating sets
respectively of $I$ and $J$. Set $I_{0} :=\sum_{i=1}^{n}Aa_{i}$
and $J_{0} :=\sum_{i=1}^{m}Ab_{i}$. Consider the exact sequence of
$R$-modules: \\
$$0 \rightarrow  Ker(u) \rightarrow  R^m  \buildrel u \over\rightarrow
J \rightarrow 0$$ where $u((c_{i},g_{i})_{1\leq i\leq m})
=\sum_{i=1}^{m}(c_{i},g_{i})(b_{i}, f_{i})$. But $Ker(u) \subseteq
(M \propto E)^m$ by \cite[Lemma 4.43, p.134]{Ro}. Hence,
\begin{eqnarray}
 \nonumber  Ker(u)&=& \{(c_{i},g_{i})_{1\leq i\leq m} \in R^m / \sum_{i=1}^{m}(c_{i},g_{i})(b_{i},f_{i}) =0\} \\
 \nonumber          &=& \{(c_{i},g_{i})_{1\leq i\leq m} \in R^m / \sum_{i=1}^{m}c_{i}b_{i} =0\} \\
\nonumber         &=& V \propto E^{m}
\end{eqnarray}
where $V :=\{(c_{i})_{1\leq i\leq m} \in A^m /
\sum_{i=1}^{m}c_{i}b_{i} =0\}$. Also, we have the exact sequence
of $R$-modules: \\
$$0 \rightarrow  V \propto E^m \rightarrow  R^m  \buildrel w \over\rightarrow
J_{0} \propto 0 \rightarrow 0$$ where $w((c_{i},g_{i})_{1\leq
i\leq m}) =\sum_{i=1}^{m}(c_{i},g_{i})(b_{i}, 0)$. Therefore, we
may assume that $J =J_{0} \propto 0$ and so $I =I_{0} \propto 0$
(since $I \subseteq J$). On the other hand, $V$ is a finitely
generated $A$-module since $V \propto E^m (=Ker(u))$ is a finitely
generated $R$-module (by the above exact sequence since $J$ is
finitely presented) and so $J_0$ is a finitely presented ideal of
$A$ by the exact sequence of $A$-modules: \\
$$0 \rightarrow  V  \rightarrow  A^m  \buildrel v \over\rightarrow
J_{0} \rightarrow 0$$ where $v((c_{i})_{1\leq i\leq m})
=\sum_{i=1}^{m}c_{i}b_{i}$. Hence, $I_0$ is a finitely presented
ideal of $A$ since $A$ is weak coherent and $I_{0} \subseteq
J_{0}$. Therefore, by the same reasoning as for $J$, we may show
that $I$ is a finitely presented ideal of $R$ and this completes
the proof of Theorem~\ref{result6}. $\Box$

\bigskip

Now, we are able to give a new class of a non-coherent weak
coherent rings. \\

\bigskip

\begin{exmp}\label{result7}
 Let $(A,M)$ be a local coherent ring with non-finitely generated maximal ideal $M$, $E$ an
 $A$-module such that $ME =0$, let $R :=A \propto E$  be the trivial ring extension of $A$
by $E$, and set $I :=0 \propto E$.  Then: \\
{\bf 1)} $R$ is a weak coherent ring by Theorem~\ref{result6}.\\
{\bf 2)} $R$ is not a coherent ring by \cite[Theorem 2.6(2)]{KM2}
since $M$ is not finitely generated. $\Box$
\end{exmp}

%%%%%%%%%%%%%%%%%%%%%%%%%%%%%%%%%%%%%%%%%%%%%%%%%%%%%%%%%%%%%%%%%%%%%%%%%%%%%%
%%%%%%%%%%%%%%%%%%%%%%%%%%%%%%%%%%%%%%%%%%%%%%%%%%%%%%%%%%%%%%%%%%%%%%%%%%%%%%%
%%%%%%%%%%%%%%%%%%%%%%%%%%%%%%%%%%%%%%%%%%%%%%%%%%%%%%%%%%%%%%%%%%%%%%%%%%%%%%
%%%%%%%%%%%%%%%%%%%%%%%%%%%%%%%%%%%%%%%%%%%%%%%%%%%%%%%%%%%%%%%%%%%%%%%%%%%%%%%
\bigskip

%%%%%%%%%%%%%%%%%%%%%%%%%%%%%%%%%%%%%%%%%%%%%%%%%%%%%%%%%%%%%
%%%%%%%%%%%%%%%%%%%%%%%%%%%%%%%%%%%%%%%%%%%%%%%%%%%%%%%%%%%%%%%

The localization of a weak coherent ring is not always a weak
coherent ring as the following Example shows:

\bigskip

%%%%%%%%%%%%%%%%%%%%%%%%%%%%%%%%%%%%%%%%%%%%%%%%%%%%%%%%%%%%%%%%%%%%%%%%%%%%%%%
%%%%%%%%%%%%%%%%%%%%%%%%%%%%%%%%%%%%%%%%%%%%%%%%%%%%%%%%%%%%%%%%%%%%%%%%%%%%%%%
\begin{exmp}\label{result8}
 Let $A :=\Z_{2} + X\R[[X]]$ a local ring with maximal ideal $M =2\Z_{2} + X\R[[X]]$,
 $E$ be an $A/M$-vector space with infinite rank and let $R :=A \propto E$ be
the trivial ring extension of $A$ by $E$. Let $S$ be the
multiplicative subset of $R$ given by $S :=\{(2,0)^{n}\ /\ n \in
\N\}$ and $S_{0}$ the multiplicative subset of $A$ given by $S_{0}
:=\{2^{n}\ /\ n \in \N\}$. Then: \\
{\bf 1)} $R$ is a weak coherent ring.\\
{\bf 2)} $S^{-1}R$ is not a weak coherent ring. $\Box$
\end{exmp}

\proof {\bf 1)} Clear by Example~\ref{result1}. \\
{\bf 2)} Since $2E \subseteq ME =0$ and $2 \in S_{0}$, then
$S_{0}^{-1}E =0$. Thus, $S^{-1}(0 \propto E) =0$ and so $S^{-1}R
=\{{(a,0)\over (s,0)}\hspace{0.2cm} / \hspace{0.2cm} a \in A$ and
$s \in S_{0}\}$ which is clearly isomorphic to a ring
$S_{0}^{-1}A$. But $S_{0}^{-1}A =S_{0}^{-1}\Z_{2} + X\R[[X]] =\Q +
X\R[[X]]$ which is not a coherent domain by \cite[Theorem
5.2.3]{Gz2}. Therefore, $S^{-1}R$ is not a weak coherent ring by
Proposition~\ref{result2}, as desired. $\Box$
%%%%%%%%%%%%%%%%%%%%%%%%%%%%%%%%%%%%%%%%%%%%%%%%%%%%%%%%%%%%%%%%%%%%%%%%%%%%%%%%%%%%%%%
%%%%%%%%%%%%%%%%%%%%%%%%%%%%%%%%%%%%%%%%%%%%%%%%%%%%%%%%%%%%%%%%%%%%%%
\bigskip

We know that a coherent ring is weak coherent and strong
$2$-coherent. The following two examples show that the class of
weakly finite conductor rings and the class of $2$-coherent rings
are not comparable. \\

\bigskip

\begin{exmp}\label{result9}
 Let $R$ be a non-coherent strong $2$-coherent domain (see for example \cite[Theorem 3.1]{M1}). Then,
 $R$ is not a weak coherent domain (by Proposition~\ref{result2} since $R$ is not a coherent domain).   $\Box$
\end{exmp}

\bigskip

\begin{exmp}\label{result10}
 Let $(A,M)$ be a local coherent domain with non-finitely generated maximal ideal $M$, and
 let $R :=A \propto (A/M)$  be the trivial ring extension of $A$
by $A/M$.  Then: \\
{\bf 1)} $R$ is a weak coherent ring by Theorem~\ref{result6}.\\
{\bf 2)} $R$ is not a strong $2$-coherent ring by \cite[Theorem
3.1]{KM1}. $\Box$
\end{exmp}

%%%%%%%%%%
%%%%%%%%%%%%%%%%%%%%%%%%%%%%%%%%%%%%%%%%%%%%%%%%%%%%%%%%%%%%%%%%%%%%%%%%%%%%%%%%%%%%%%%
%%%%%%%%%%%%%%%%%%%%%%%%%%%%%%%%%%%%%%%%%%%%%%%%%%%%%%%%%%%%%%%%%%%%%%
%%%%%%%%%%%%%%%%%%%%%%%%%%%%%%%%%%%%%%%%%%%%%%%%%%%%%%%%%%%%%%%%%%%%%%%%%%%%%%%%%%%%%%%
\bigskip

 Next, we study the transfer of weak coherent property to direct products. \\

\bigskip

\begin{prop}\label{result11}
 Let $(R_{i})_{i=1,\ldots ,n}$ be a family of rings. Then, $\prod _{i=1}^{n}R_{i}$ is a weak coherent ring
 if and only if $R_i$ is a weak coherent ring for each $i =1, \ldots ,n$. $\Box$
\end{prop}

\bigskip

We need the following Lemma before proving Proposition~\ref{result11}. \\

\bigskip

%%%%%%%%
\begin{lem} (\cite[Lemma 2.5(1)]{M1})\label{result12}  Let $(R_{i})_{i=1,2}$ be a family of
rings and $E_i$ an $R_i$-module for $i =1,2$. Then  $E_{1} \prod
E_{2}$ is a finitely generated (resp., finitely presented) $R_{1}
\prod R_{2}$-module if and only if
$E_i$ is a finitely generated (resp., finitely presented) $R_i$-module for $i =1,2$. \\

 \end{lem}

\bigskip

 %%%%%%%%%%%%%%%%%%%%%%%%%%%%%%%%%%%%%%%%%%%
{\bf Proof of Proposition 2.11.} By induction on $n$, it suffices
to prove the assertion for $n =2$. Since an ideal of $R_{1} \prod
R_{2}$ is of the form $I_{1} \prod I_{2}$, where $I_i$ is an ideal
of $R_i$ for $i =1, 2$, the conclusion follows easily from
Lemma~\ref{result11}.  \\
%%%%%%%%%%
\bigskip

Now, we are able to give a new class of a non-coherent weak
coherent rings. \\

\bigskip

\begin{exmp}\label{result13}
 Let $R_1$ be a non-coherent weak coherent ring, $R_2$ be a coherent ring, and $R =R_{1} \prod R_{2}$.  Then: \\
{\bf 1)} $R$ is a weak coherent ring by Theorem~\ref{result10}.\\
{\bf 2)} $R$ is not a coherent ring by \cite[Theorem 2.4.3]{Gz2}.
$\Box$
\end{exmp}
%%%%%%%%%%%%%%%%%%%%%%%%%%%%%%%%%%%%%%%%%%%%%%%%%%%%%%%%%%%%%%%%%%%%%%%%%%%%%%%%%%%%%%%
%%%%%%%%%%%%%%%%%%%%%%%%%%%%%%%%%%%%%%%%%%%%%%%%%%%%%%%%%%%%%%%%%%%%%%
%%%%%%%%%%%%%%%%%%%%%%%%%%%%%%%%%%%%%%%%%%%%%%%%%%%%%%%%%%%%%%%%%%%%%%%%%%%%%%%%%%%%%%%

\end{section}

%%%%%%%%%%%%%%%%%%%%%%%%%%%%%%%%%%%%%%%%%%%%%%%%%%%%%%%%%%%%
%%%%%%%%%%%%%%%%%%%%%%%%%%%%%%%%%%%%%%%%%%%%%%%%%%%%%%%%%
%%%%%%%%%%%%%%%%%%%%%%%%%%%%%%%%%%%%%%%%%%%%%%%%%%%%%%%%%
%%%REFERENCES%%%%%%%%%%%%%%%%%%%%%%%%%%%%%%%%%%%%%%%%%%%%
%%%%%%%%%%%%%%%%%%%%%%%%%%%%%%%%%%%%%%%%%%%%%%%%%%%%%%%%

%%%%%%%%%%%%%%%%%%%%%%%%%%%%%%%%%%%%%%%%%%%%%%%%%%%%

%%%%%%%%%%%%%%%%%%%%%%%%%%%%%%%%%%%%%%%%%%%%%%%%%%%%

\bigskip\bigskip

%%%%%%%%%%%%%%%%%%%%%%%%%%%%%%%%%%%%%%%%%%%%%%%%%%%%%%%%

\begin{thebibliography}{999}\addcontentsline{toc}{section}{\protect\numberline{}{Bibliography}}





\bibitem{Bo} N. Bourbaki, {\em Commutative algebra, Chapters 1-7}, Springer,
Berlin, 1998.

\bibitem{C} D. Costa, {\em Parameterizing families of non-Noetherian rings}, Comm.
Algebra 22 (1994), 3997--4011.

\bibitem{Gi} R. Gilmer, {\em Multiplicative ideal theory}, Pure and Applied
Mathematics, No. 12. Marcel Dekker, Inc., New York, 1972.

\bibitem{Gz1}  S. Glaz, {\em Finite conductor rings}, Proc. Amer. Math. Soc. 129 (2000), 2833--2843.

\bibitem{Gz2}  S. Glaz, {\em Commutative coherent rings}, Lecture Notes in
Mathematics, 1371, Springer-Verlag, Berlin, 1989.

\bibitem{Gz3}  S. Glaz, {\em Controlling the Zero-Divisors of a Commutative Ring},  Lecture Notes in Pure and Appl. Math.,
Dekker, 231 (2003), 191--212.

\bibitem{H}  J.A. Huckaba, {\em Commutative rings with zero divisors}, Marcel
Dekker, New York-Basel, 1988.

\bibitem{KM1} S. Kabbaj, N. Mahdou, {\em Trivial extensions of local rings and a
conjecture of Costa}, Lecture Notes in Pure and Appl. Math.,
Dekker, 231 (2003), 301--311.

\bibitem{KM2} S. Kabbaj and N. Mahdou, {\em Trivial extensions defined by coherent-like conditions}, Comm. Algebra 32(10) (2004), 3937--3953.

\bibitem{M1} N. Mahdou, {\em On Costa's conjecture}, Comm. Algebra 29
(2001), 2775--2785.

\bibitem{M2} N. Mahdou, {\em On 2-Von Neumann regular rings}, Comm. Algebra 33 (10) (2005) 3489--3496.

\bibitem{Ro} J. J. Rotman, {\em An Introduction to Homological Algebra}, Academic Press, New York, 1979.


    \end{thebibliography}
\end{document}